\documentclass[11pt,letterpaper]{article}
\usepackage{amsfonts, amsmath, amssymb, amscd, amsthm, color, graphicx, mathabx, mathrsfs, wasysym, setspace, mdwlist, calc, float, mathtools, tikz-cd}

\usepackage{tocloft}
\usepackage{xcolor}

\usepackage{hyperref}
\numberwithin{equation}{section}

\hypersetup{colorlinks,
    linkcolor={red!50!black},
    citecolor={blue!80!black},
    urlcolor={blue!80!black}}
\usepackage{float}
\newcommand{\Ga}{\Gamma}
\newcommand{\La}{\Lambda}

\newcommand{\ca}{\curvearrowright}

\newcommand{\A}{\mathcal A}
\newcommand{\B}{\mathcal B}
\newcommand{\C}{\mathcal C}

\newcommand{\G}{\mathcal G}

\renewcommand{\L}{\mathcal L}
\newcommand{\M}{\mathcal M}
\newcommand{\N}{\mathcal N}
\renewcommand{\P}{\mathcal P}
\newcommand{\Q}{\mathcal Q}

\newcommand{\V}{\mathcal V}

\newcommand{\sD}{\mathscr D}
\newcommand{\sE}{\mathscr E}
\newcommand{\sF}{\mathscr F}
\newcommand{\sG}{\mathscr G}
\newcommand{\sH}{\mathscr H}

\newcommand{\sK}{\mathscr K}

\newcommand{\sN}{\mathscr N}

\newcommand{\sS}{\mathscr S}

\newcommand{\sU}{\mathscr U}
\newcommand{\sV}{\mathscr V}
\newcommand{\sW}{\mathscr W}

\newtheorem{thm}{Theorem}[section]
\newtheorem*{thm*}{Theorem}
\newtheorem{cor}[thm]{Corollary}
\newtheorem{lem}[thm]{Lemma}

\theoremstyle{definition}

\theoremstyle{remark}

\begin{document}

\title{Cartan Subalgebras in von Neumann Algebras Associated with Graph Product Groups}
\date{\today}
\author{Ionu\c t Chifan and Srivatsav Kunnawalkam Elayavalli}




\maketitle

\begin{abstract} In this short note we classify the Cartan subalgebras in all von Neumann algebras associated with graph product groups and their free ergodic measure preserving actions on probability spaces. 
    
\end{abstract}

\section{Introduction}
Murray and von Neumann found a natural way to associate a crossed-product von Neumann algebra, denoted by $L^\infty(X)\rtimes \Gamma$, to any measure preserving action $\Gamma \ca X$ of a countable group $\Gamma$ on a probability space $X$ \cite{MvN36}. When $X$ consists of a singleton this amounts to the group von Neumann algebra $\L(\Gamma)$ \cite{MvN43}. The most interesting cases for study is when  $L^\infty(X)\rtimes \Gamma$ and  $\L(\Gamma)$ are II$_1$ factors which correspond to the situations when $\Gamma \ca X$ is a free, ergodic action and $\Gamma$ is an icc group, respectively.   

An important open problem regarding the structure of these algebras is to classify their Cartan subalgebras\footnote{$\A\subset \M$ is a Cartan subalgebra if it is a maximal abelian von Neumann subalgebra and its normalizer $\mathcal N_{\M}(\A)=\{ u\in \mathcal U(\M) \,:\, u\A u^*=\A\}$ generates the entire $\M$, as a von Neumann algebra.} for various natural choices of groups/actions. Specifically, the main interest is to identify situations when these algebras have no Cartan subalgebras or have a unique Cartan subalgebra, up to unitary conjugacy. The first progress on this problem was made in \cite{Po81} where Popa showed the (nonseparable) factors arising from the free group with uncountably many generators do not have Cartan subalgebras.  Using his influential free entropy dimension theory Voiculescu was able to show that the free group factors with countably many generators do not have Cartan subalgebras \cite{Vo96}. Over the last $15$ years we have witnessed an unprecedented progress towards this problem through the emergence of Popa deformation/rigidity theory. Using this powerful conceptual framework Ozawa and Popa discovered that the free group factors are strongly solid (in particular have no Cartan subalgebras) and all crossed products $L^\infty(X)\rtimes \mathbb F_n$ arising from free ergodic profinite actions $\mathbb F_n \ca X$ with $n\geq 2$ have unique Cartan subalgebra \cite{OP07}. The breakthrough technology developed in this work laid out the foundations for a plethora of subsequent results regarding the study of Cartan subalgebras in II$_1$ factors, \cite{OP08,CS11,PV11,PV12,Io12,Va13,CIK13}---just to name a few. In particular, the developments from \cite{PV11} due to Popa and Vaes had a profound impact in the study of normalizers of amenable von Neumann algebras and the entire program of classification of von Neumann algebra arising from groups and their actions on probability spaces.

In this paper we will investigate the aforementioned problems when $\Gamma$ is a graph product of non trivial groups. These are natural generalizations of free products and direct products and were introduced by E.\ Green in her PhD thesis \cite{Gr90}. Their study has became a trendy subject over the last decades in topology and geometric group theory with many important results emerging in recent years. Important and well studied examples of such groups are the right angled Artin and Coxeter groups.  
General graph product groups have been considered in the von Neumann algebras setting for the first time in \cite{CF14} and several structural results have been established for some of them \cite{CF14,Ca16,CdSS18,DK-E21}. In \cite{Ca16} Caspers was able to isolate fairly general conditions on a graph that are sufficient to insure the corresponding Coxeter group $\Gamma$ gives rise to a von Neumann algebra $\L(\Gamma)$ with no Cartan subalgebras. The scope of this note is to broaden this study by providing necessary/sufficient conditions on the graph product data $\Gamma=\mathcal G\{\Gamma_v\}$ that will allow the classification of the Cartan subalgebras in $\L(\Gamma)$ and $L^\infty(X)\rtimes \Gamma$.

To properly introduce our results we briefly recall the graph product group construction. Let $\sG=(\sV,\sE)$ be a finite, simple (no self loops or multiple edges) graph, where $\sV$ and $\sE$  denote its vertices and edges sets, respectively. Let $\{\Ga_v\}_{v\in\sV}$ be a family of groups called vertex groups. The graph product group associated with this data, denoted by  $\sG\{\Ga_v\}$, is the group generated  by $\Ga_v$, $v\in \sV$ subject to the relations of the groups $\Gamma_v$ along with the relations $[\Ga_u, \Ga_v] = 1$, whenever $(u,v)\in \sE$. Developing an approach that combines deep results in Popa's deformation/rigidity theory with various geometric/algebraic aspects of graph product groups we are able to classify the Cartan subalgebras for large classes of von Neumann algebras associated with graph product groups. Specifically, using the classification of normalizers of amenable subalgebras in amalgamated free product von Neumann algebras due to Ioana \cite{Io12} (see also \cite{Va13}) in cooperation with the transitivity of intertwining virtual Cartan subalgebras from \cite{Io11,CIK13} and the structure of the parabolic subgroups of a graph product groups \cite{AM10} we have the following:

\begin{thm}\label{nocartan1}Let $\Ga= \sG\{\Ga_v\}$ be an icc graph product of groups where $\sG$ does not admit a de Rham join decomposition of the form  $\sG=\sG_1\circ \sG_2\circ \cdots \circ \sG_k$ so that for all $2\leq i\leq k$ we have $\sG_i=(\{v_i,w_i\}, \emptyset)$ and $\Ga_{v_i}=\Ga_{w_i}\cong \mathbb Z_2$. Then $\L(\Ga)$ has no Cartan subalgebra.

\end{thm}

In particular, when the graph product does not contain elements of order two (i.e.\ non-Coxeter type) we get a result which complements the results from \cite{Ca16}.  

\begin{cor}\label{nocartan2}Let $\sG=(\sV,\sE)$ be any finite simple graph such that there exist $v,w\in \sV$ satisfying $(v,w)\notin \sE$. Assume that  $\Ga=\sG(\Ga_v)$ is any icc graph product group where $|\Ga_v|\geq 3$ for all $v\in \sV$. Then $\L(\Ga)$ has no Cartan subalgebra.

\end{cor}

In their groundbreaking work \cite{PV11}, Popa and Vaes discovered the first examples of group $\Gamma$ that are $\mathcal C$-rigid.  This means that any free ergodic pmp action $\Gamma\ca X$ gives rise to a von Neumann algebra $L^\infty(X)\rtimes \Gamma$ whose ``functions'' subalgebra  $L^\infty(X)$ is the unique Cartan subalgebra, up to unitary conjugacy. Remarkably, their examples include all direct products of free groups $\Gamma=\mathbb F_{n_1}\times \mathbb F_{n_2}\times \cdots \times\mathbb F_{n_k}$, for all $n_i\geq 2$ and $k\geq 1$ as well as other examples of right-angled Artin groups. In this paper we are able to provide additional examples  of $\mathcal C$-rigid general graph products groups thus complementing these results.

\begin{thm}\label{uniquecartan1} Let $\Ga= \sG\{\Ga_v\}$ be a graph product of groups so that the graph $\sG$ does not admit a nontrivial join decomposition $\sG= \sH\circ \sK $ where either $\sH=(\{v\},\emptyset)$, or, $\sH=(\{v,w\}, \emptyset)$ with $\Ga_v\cong\Ga_w\cong \mathbb Z_2$. Let $\Ga\ca X $ be any free ergodic pmp action  and let $\M = L^\infty(X)\rtimes \Ga$ be the corresponding group measure space von Neumann algebra. Then for any Cartan subalgebra  $\A\subset \M $ there is a unitary $u\in \M$ so that $\A =u L^\infty(X) u^*$.  
\end{thm}

While the graph product conditions highlighted in the statement appear little more restrictive than the ones presented in Theorem \ref{nocartan1}, they are  in fact optimal. The reader may consult to the remarks at the end of the paper for several counterexamples in this direction which build on  \cite{OP07}. Hence, the results of this paper offer a complete solution to Question 1.9 raised by the second author and C. Ding in \cite{DK-E21}.

\section{Preliminaries}

\subsection{Graph products of groups}

In this subsection we briefly recall the notion of graph product of groups introduced by E. Green  \cite{Gr90} also highlighting some of its features that are relevant to this note. Let $\sG=(\sV,\sE)$ be a finite simple graph, where $\sV$ and $\sE$  denote its vertices and edges sets, respectively. Let $\{\Ga_v\}_{v\in\sV}$ be a family of groups called vertex groups. The graph product group associated with this data, denoted by $\sG\{\Ga_v,v \in \sV\}$ or simply $\sG\{\Ga_v\}$, is the group generated  by $\Ga_v$, $v\in \sV$  with the only relations that $[\Ga_u, \Ga_v] = 1$, whenever $(u,v)\in \sE$. 
Given any subset $\sU\subset \sV$, the subgroup $\Ga_\sU =\langle \Ga_u \,:\,u\in \sU\rangle $ of $\sG\{\Ga_v,v\in \sV\}$ is called a full subgroup. This can be identified  with the graph product $\sG_\sU\{\Ga_u,u \in \sU\}$ corresponding to the subgraph $\sG_\sU$ of $\sG$, spanned by the vertices of $\sU$. 
For every $v \in \sV$ we denote by ${\rm lk}(v)$ the subset of vertices $w\neq v$ so that $(w,v)\in \sE$. Similarly, for every $\sU 
\subseteq \sV$ we denote by ${\rm lk}(\sU) = \cap _{u\in \sU}{\rm lk}(u)$. Also we make the convention that ${\rm lk}(\emptyset) = \sV$. Notice that $\sU \cap  {\rm lk}(\sU) = \emptyset$.

Graph product groups admit many amalgamated free product decompositions. One such decomposition which is essential on deriving our main results, involves full subgroups factors in   \cite[Lemma 3.20]{Gr90} as follows. For any $w \in \sV$ we have
$$\sG\{\Ga_v\} = \Ga_{\sV\setminus \{w\}} \ast_{ \Ga_{\rm lk}(w)} \Ga_{{\rm st}(w)},$$
where ${\rm st} (w) = \{w\} \cup {\rm lk} (w)$. Notice that always $\Ga_{{\rm lk}(w)}\lneqq\Ga_{{\rm st}(w)}$ but it could be the case that $\Ga_{{\rm lk}(w)}=\Ga_{\sV\setminus \{w\}} $, when $\sV={\rm st}(w)$. In this case the amalgam decomposition is degenerate.

If a graph $(\sV,\sE)=\sG$ has two proper subgraphs $(\sW,\sF)=\sH,(\sU,\sD)=\sK \subset \sG$ whose vertices form a partition $\sV=\sW\sqcup \sU$ and $(v,w)\in \sE$ for all $v
\in \sW$ and $w\in \sU$ then we say that $\sG$ admits a non-trivial join product decomposition $\sG = \sH\circ \sK$. A graph $\sG=(\sV,\sE)$ is called irreducible if it does not admit any non-trivial join product decomposition. It is well-known that any graph $\sG$ admits a (de Rham) join product decomposition $\sG= \sG_0\circ\sG_1 \circ \cdots \circ \sG_k$ for $k\geq 0$, where $\sG_0$ is a clique (i.e.\ a complete graph) and $\sG_i\subseteq \sG$, for all $i\geq 1$ are irreducible subgraphs with at least two vertices.  

\subsection {Popa's intertwining-by-bimodules techniques} We next recall from  \cite [Theorem 2.1, Corollary 2.3]{Po03} Popa's {\it intertwining-by-bimodules} technique.
Let $\Q\subset \M$ be a  von Neumann subalgebra. The  basic construction $\langle \M,e_\Q\rangle$ is defined as the von Neumann subalgebra of $\mathbb B(L^2(\M))$ generated by $\M$ and the orthogonal projection $e_\Q$ from $L^2(\M)$ onto $L^2(\Q)$. There is a semi-finite faithful trace on $\langle \M,e_\Q\rangle$ given by $\text{Tr}(xe_\Q y)=\tau(xy)$, for every $x,y\in \M$. We denote by $L^2(\langle \M,e_\Q\rangle)$ the associated Hilbert space and endow it with the natural $\M$-bimodule structure.

\begin{thm}[\cite{Po03}]\label{corner} Let $(\M,\tau)$ be a tracial von Neumann algebra and $\P\subset p\M p,\Q\subset \M$ be von Neumann subalgebras. 
Then the following  are equivalent:

\begin{enumerate}
\item There exist projections $p_0\in \P, q_0\in \Q$, a $*$-homomorphism $\theta:p_0\P p_0\rightarrow q_0\Q q_0$  and a non-zero partial isometry $v\in q_0\M p_0$ such that $\theta(x)v=vx$, for all $x\in p_0\P p_0$.

\item There exists a non-zero element $a\in \P'\cap p\langle\M,e_\Q\rangle p$ such that $a\geq 0$ and $\text{Tr}(a)<\infty$.

\item There is no sequence $u_n\in\mathcal U(\P)$ satisfying $\|E_\Q(x^*u_ny)\|_2\rightarrow 0$, for all $x,y\in p\M$.
\end{enumerate}
If one of these equivalent conditions holds,  we write $\P\prec_{\M}\Q$, and say that \emph{a corner of $\P$ embeds into $\Q$ inside $\M$.}
Moreover, if $\P p'\prec_{\M}\Q$ for all projections $0\neq p'\in \P'\cap p\M p$, then we write $\P\prec^{s}_{\M}\Q$.\end{thm}
Next we recall a few intertwining results that will be used in the main proofs. The first is a straightforward generalization of \cite[Lemma 2.]{CI17}, its proof being left to the reader.
\begin{lem}[\cite{CI17}]\label{groupintertwiner} Let $\La, \Sigma\leqslant \Ga$ be groups, let $\Ga \ca \P$ be a  trace preserving action on a von Neumann algebra and let $\M =\P\rtimes \Ga$ be the corresponding crossed product von Neumann algebra. Then $\P\rtimes \La \prec_\M \P \rtimes \Sigma$ if an only if there is $g\in \Ga$ such that $[\La: \La \cap g \Sigma g^{-1}]<\infty$. 

\end{lem}

The second is Popa's  conjugacy criterion for Cartan subalgebras \cite[Theorem A.1]{Po01}.

\begin{thm}[\cite{Po01}]\label{Po01} Let $\M$ be a II$_1$ factor and $\P,\Q\subset \M$ be Cartan subalgebras. If $\P\prec_\M\Q$, then there is $u\in\sU(\M)$ such that $u\P u^*=\Q$. 
\end{thm}

\subsection{Relative amenability}

Recall that the tracial von Neumann algebra $(\M,\tau)$  is called {\it amenable} if there exists a sequence $\xi_n\in L^2(\M)\otimes L^2(\M)$ such that $\langle x\xi_n,\xi_n\rangle\rightarrow\tau(x)$ and $\|x\xi_n-\xi_nx\|_2\rightarrow 0$, for every $x\in \M$. By Connes' celebrated classification of injective factors \cite{Co76}, $\M$ is amenable if and only if it is approximately finite dimensional.

Let $\P\subset p\M p$ and $\Q \subseteq \M$ be a von Neumann subalgebras.
Following \cite[Section 2.2]{OP07} we say that $\P$ is  {\it amenable relative to $\Q$ inside $\M$} if there exists a sequence $\xi_n\in L^2(\langle \M,e_\Q\rangle)$ such that $\langle x\xi_n,\xi_n\rangle\rightarrow\tau(x)$, for every $x\in p\M p$, and $\|y\xi_n-\xi_ny\|_2\rightarrow 0$, for every $y\in \P$. When $\P=\M$ one simply says that $\M$ is \emph{co-amenable relative to} $\Q$, 
 paralleling Eymard's notion of co-amenability in group theory---a subgroup $H$ of $G$ is co-amenable if there exists a left $G$ invariant state on $\ell^{\infty}(G/H)$.  


The following lemma is well known, but for convenience of the reader we include a short proof below:

\begin{lem}\label{coamensubgr} Let $\Ga=\Ga_1\ast_\Sigma\Ga_2$ be a non-degenerate amalgamated free product such that $\Ga$ is co-amenable relative to $\Sigma$. Then $[\Ga_i:\Sigma]= 2$ for all $i=1,2$.

\end{lem}
\begin{proof}
   Suppose that $[\Gamma_1:\Sigma]\geq3$, and suppose $\varphi$ is a non trivial left $\Gamma$-invariant state on $\ell^{\infty}(\Gamma/\Sigma)$. Let $S_i$ denote the subset of elements in $\Gamma/\Sigma$ whose normal form decomposition begins with $\Gamma_i/\Sigma$. Clearly, we have $\phi(\chi_{S_1}+\chi_{S_2})=1$. Moreover, since $[\Gamma:\Gamma_1]\geq3$, there exists $g_1\neq g_2\in \Gamma/\Gamma_1$. Let $\Gamma_2/\Sigma\ni h\neq 1 $. Hence $ \varphi(g_1\chi_{S_2} +g_2\chi_{S_2})\leq \varphi(h\chi_{S_1})\leq  \varphi(\chi_{S_2})$, giving $\varphi(\chi_{S_2})=0$. Similarly $\varphi(\chi_{S_1})=0$ which is a contradiction.
\end{proof}

\section{Proofs of the main results}

The main technical result of the paper is the following intertwining result for \emph{virtual} Cartan subalgebras in crossed products by graph product groups. All other results are derived from it. Its proof relies on the classification of normalizers of amenable subalgebras in amalgamated free products due to Ioana \cite{Io12} and its more general version by Vaes \cite{Va13} together with a technique on transitivity of intertwining virtual Cartan subalgebras from  \cite{Io11,CIK13}.

\begin{thm}\label{controlmasa1} Let $\Ga =\sG \{\Ga_v\}$ be a graph product of groups and let $\G \ca \P$ be a trace preserving action on a von Neumann algebra.  
Denote by $\M= \P\rtimes \Ga$ the corresponding crossed product von Neumann algebra and let  $p\in \M$ be a projection. Assume that $\A\subset p\M p$ is a masa such that its normalizer $\sN_{p\M p}(\A)''\subseteq p \M p $ has finite index. Thus one can find a subgraph $\sG_0 \subseteq \sG$ whose de Rham join product decomposition $\sG_0=\sG_1\circ \sG_2\circ \cdots \circ \sG_k$ satisfies $\sG_i =(\{v_i,w_i\},\emptyset)$ with $\Ga_{v_i}=\Ga_{w_i}\cong \mathbb Z_2$ for all $2\leq i\leq k$ such that $\A \prec_\M \P \rtimes \Ga_{\sG_0}$.\end{thm}

\begin{proof}
Let $(\sV_0,\sE_0)=\sG_0\subseteq \sG$ be a subgraph such that $|\sV_0|$ is minimal with the property that $\A \prec_\M \P\rtimes \Ga_{\sG_0}=: \N$. Let $\sG_0=\sG_1\circ \sG_2\circ \cdots \circ \sG_k$ be its de Rham decomposition where $\sG_1$ is a clique of $\sG_0$. Next we show all $\sG_i$ with $i\geq 2$ have the properties listed in the statement. 

Since $\A$ is a masa, using  \cite[Proposition 3.6]{CIK13} there are projections $0\neq a\in \A,0\neq e\in \N$, a masa $\B\subset e\N e$, a projection $0\neq e'\in \B'\cap e\M e$ and a unitary $u\in \M $ such that \begin{equation}\label{findeximage}\begin{split}& \text{the inclusion }\Q:=\sN_{e \N e}(\B)''\subseteq e \N e \text{ has finite index,}\\ & \text{the support satisfies }s(E_\N(e'))=e \text{ and } u\A a u^*=\B e'.\end{split}\end{equation} 

Fix $i\geq 2$ and let $v\in \V_i$ such that $\Ga_{\sG_0}$ admits a non-canonical decomposition $\Ga_{\sG_0}=\Ga_{\sG_0\setminus \{v\}}\ast_{\Ga_{\text{lk}(v)}} \Ga_{\text {st}(v)}$; in particular we have that $|\text{lk}(v)|\leq|\sV_0|-2$. Using \cite{Va13, Io12} one of the following cases must hold:
\begin{enumerate}
    \item [i)] $\B\prec_{\N} \P \rtimes \Gamma_{\text{lk}(v)}$;
    \item [ii)] $\Q \prec_{\N} \P \rtimes \Ga_{\sG_0\setminus \{v\}})$ or $\P\rtimes \Ga_{\text{st}(v)}$;
    \item [iii)] $\Q$ is amenable relative to $\P \rtimes \Ga_{\text{lk}(v)}$ inside  $\N$.
\end{enumerate}

Assume i). Thus one can find a partial isometry $0\neq w\in \N$ with $w^*w\in \B$ such that $w\B \subseteq (\P\rtimes \G_{\text{lk}(v)})w$. Now we argue that $we'\neq 0$. If $0=we'$ then $0=w^*we'$ and since $w\in \N$ we get $0=w^*wE_\N(e')$.  This however further implies that $0= w^*w s(E_\N(e'))= w^*w e=w^*w$, which is contradiction. Therefore we have $w e'\B =w \B e'\subseteq (\P\rtimes \G_{\text{lk}(v)})we'$ and using the second part in \eqref{findeximage} we get $w u a\A \subseteq (\P\rtimes \G_{\text{lk}(v)})we'u$. In particular, $\A \prec_\M \P\rtimes \Gamma_{\text{lk}(v)}$ and since $|\text{lk}(v)|\leq |\sV_0|-1$, this contradicts the minimality of $|\sV_0|$. Thus i) cannot hold.

Now assume $\Q \prec_{\N} \P\rtimes \Ga_{\sG_0\setminus \{v\}}$. Using \cite[Lemma 2.4(3)]{DHI16} there is a projection $0\neq q\in \Q' \cap e\N e$ such that $\Q q \prec^{\text{s}}_{\N}\P\rtimes \Ga_{\sG_0\setminus \{v\}} $.  Part one in  \eqref{findeximage} implies that $\Q q\subseteq  q\N q$ has finite index. Then combining this with \cite[Lemma 2.2]{CIK13} and \cite[Lemma 2.4(1)]{DHI16} we conclude that $\P\rtimes \Ga_{\sG_0}\prec_{\N} \P\rtimes \Ga_{\sG_0\setminus \{v\}}$. By Lemma \ref{groupintertwiner} it follows that $\Ga_{\sG_0\setminus \{v\}}\leqslant \Ga_{\sG_0}$ has finite index, which is a contradiction. To see this pick $t\in \sV_0\setminus \{v\}$ such that $(t, v)\notin \sE$. Now for any $g\in \Ga_t\setminus \{1\}$ and $h\in \Ga_v\setminus\{1\}$ one can see that $(gh)^k \Ga_{\sV_0\setminus \{v\}}$ for all $k\geq 2$ are distinct co-sets in $\Ga_{\sV_0}$.     If $\Q \prec_{\N} \P\rtimes \Ga_{\text{st}(v)}$ then in a similar manner one gets that $\Ga_{\text{st}(v)}\leqslant \Ga_{\sG_0}$ has finite index which is again contradiction. Altogether these show  ii) cannot hold either.

Assume iii). Since $\Q\subseteq e\N e$ has finite index we have $e\N e$ is amenable relative to $\Q$ inside  $\N$ and thus by \cite[Proposition 2.3]{OP07} we get that $e\N e$ is amenable relative to $\L(\Ga_{\text{lk}(v)})$. This implies that $\Ga_{\sG_0}$ is co-amenable relative to $\Ga_{\text{lk}(v)}$. Thus by Lemma \ref{coamensubgr} we have $\sG_0\setminus \text{st}(v) =\{t\}$, $\Ga_t,\Ga_v\cong \mathbb Z_2$ and also $\text{lk}(v)=\text{lk}_{\sG_0}(t)$. Thus, a fortiori we have $\sG_i =(\{v,t\}, \emptyset)$ which gives the desired conclusion. \end{proof}

\begin{thm}\label{controlmasa2} Let $\Ga =\sG \{\Ga_v\}$ be a graph product of groups and let $\G \ca \P$ be a trace preserving action on a von Neumann algebra.  
Assume the corresponding crossed product  $\M= \P\rtimes \Ga$ is a II$_1$ factor and let  $p\in \M$ be a projection. Assume that $\A\subset p\M p$ is a Cartan subalgebra. Then one can find a join decomposition $\sG= \sG_0 \circ \sK$ where $\sG_0$ has de Rham join product decomposition $\sG_0=\sG_1\circ \sG_2\circ \cdots \circ \sG_k$ so that for all $2\leq i\leq k$ we have $\sG_i=(\{v_i,w_i\}, \emptyset)$ and $\Ga_{v_i}=\Ga_{w_i}\cong \mathbb Z_2$ which satisfies $\A \prec_\M \P \rtimes \Ga_{\sG_0}$.

 \end{thm}

\begin{proof} Consider $\sG = \sH_0\circ \sH_1 \circ \cdots \circ\sH_l$ be the de Rham join product decomposition of $\sG$ and let $\sH_i=(\sV_i, \sE_i)$. Form Theorem  \ref{controlmasa1}  there is a subgraph $\sG_0\subseteq \sG$ with de Rham join product decomposition $\sG_0=\sG_1\circ \sG_2\circ \cdots \circ \sG_k$ so that for all $i\geq 2$ we have $\sG_i=(\{v_i,w_i\}, \emptyset)$ and $\Ga_{v_i}=\Ga_{w_i}\cong \mathbb Z_2$ which satisfies $\A \prec_\M \P \rtimes \Ga_{\sG_0}$. Denote by $\sW_i=\{v_i,w_i\}$. Using \cite[Lemma 2.4(3)]{DHI16} we have \begin{equation}\label{int2}
   A \prec^s_\M \P \rtimes \Ga_{\sG_0}. 
\end{equation}
Moreover, we can assume $\sG_0$ is a subgraph with all the properties listed above and a minimal number of vertices $|\sV_0|$. In particular this implies that $\sG_0$ cannot be represented as a link $\sG_0={\rm lk}(w)$ of a vertex such that $\Ga_w$ is a finite group.

Next let $I_1\sqcup I_2 \sqcup \cdots \sqcup I_s= \{1,\ldots, k\}$ be a partition and let $\{i_1,i_2, \ldots, i_t\}\subseteq \{0,\ldots, l\}$ be a subset such that   $\circ_{j\in I_r}\sG_j\subseteq \sH_{i_r}$  is a subgraph for $1\leq r\leq t$. Here we have $s\leq k$ and $t\leq l+1$. Notice that if $i_1=0$ then necessarily $I_1=\{1\}$. 

Let $i_r\neq 0$ and assume $ \cup_{j\in I_r} \sW_{j}  \subsetneq \sV_{i_r} $. Since $\sG_{i_r}$ is irreducible there is $v\in\sV_{i_r}\setminus  \cup_{j\in I_r} \sW_{j}$ such that ${\rm lk}_{\sG_{i_r}}(v)\cap \circ_{j\in I_r} \sG_j\subsetneq \circ_{j\in I_r} \sG_j$. In particular, we also have ${\rm lk}(v) \cap \sG_0\subsetneq \sG_0$. Now we note that $\Ga_{\sG}$ admits a non-canonical decomposition $\Ga_{\sG}=\Ga_{\sG\setminus \{v\}}\ast_{\Ga_{\text{lk}(v)}} \Ga_{\text {st}(v)}$. Using \cite{Va13, Io12} and \cite[Lemma 2.4(3)]{DHI16} one of the following cases must hold:
\begin{enumerate}
    \item [i)] $\A\prec^s_{\M} \P \rtimes \G_{\text{lk}(v)}$;
    \item [ii)] $p\M p\prec^s_{\M} \P \rtimes \Ga_{\sG\setminus \{v\}})$ or $\P\rtimes \Ga_{\text{st}(v)}$;
    \item [iii)] $p\M p$ is co-amenable relative to $\P \rtimes \Ga_{\text{lk}(v)}$.
\end{enumerate}
Assume i) holds. Combining it with \eqref{int2} and using \cite[Lemma 2.6]{DHI16} there is $h\in\Ga$  such that $\A \prec^s \P\rtimes \left(\Ga_{{\rm lk}(v)} \cap h \Ga_{\sG_0} h^{-1}\right )$. By \cite[Proposition 3.4]{AM10} there exists $\sS\subset {\rm lk}(v)\cap \sG_0\subsetneq \sG_0$ such that $k\Ga_\sS k^{-1}=\Ga_{{\rm lk}(v)} \cap h \Ga_{\sG_0} h^{-1}$ for some $k\in \Ga$  and hence $\A \prec^s \P\rtimes \Ga_\sS$. However the graph $\sS$ has all the properties of $\sG_0$ but but it has a smaller number of vertices thus contradicting the minimality of $|\sV_0|$. So i) cannot hold. 

Assume ii) holds. Then by Lemma \ref{groupintertwiner} we have either $[\Ga: \Ga_{\sG\setminus \{v\}}]<\infty$ or $[\Ga: \Ga_{{\rm st}(v)}]<\infty $ which are impossible.

Assume iii) holds. Therefore $\Ga$ is amenable relative to $\Ga_{\rm lk(v)}$ and by Lemma \ref{coamensubgr} we have $\sG\setminus \text{st}(v) =\{t\}\in \circ_{j\in I_r} \sH_j$, $\Ga_t,\Ga_v\cong \mathbb Z_2$ and thus $\sG= \{v,t\}\circ \text{lk}\{v,t\}$. We also have that $\{v,t\}\subseteq \sV_{i_r}$ and since $\sG_{i_r}$ is irreducible we conclude that $\sG_{i_r}=(\{v,t\}, \emptyset)$. But this implies in particular that $(\{t\},\emptyset)= \circ_{j\in I_r} \sH_j$ ans since $\Ga_t\cong \mathbb Z_2$ this again contradicts the minimality of $|\sV_0|$. So this case is impossible as well.

 In conclusion, for all $i_r\neq 0$ we have   $\circ_{j\in I_r} \sH_j =\sG_{i_r}$ which gives the desired conclusion. \end{proof}
\vskip 0.03in
\noindent \emph{Proof of Theorem \ref{nocartan1}.} Assume by contradiction $\A \subset \M =\L(\Ga)$ is a Cartan subalgebra. From Theorem \ref{controlmasa2} there is a decomposition $\sG= \sG_0 \circ {\rm lk}(\sG_0)$ with de Rham join product decomposition $\sG_0=\sG_1\circ \sG_2\circ \cdots \circ \sG_k$ so that for all $2\leq i\leq k$ we have $\sG_i=(\sV_i =\{v_i,w_i\}, \emptyset)$ and $\Ga_{v_i}=\Ga_{w_i}\cong \mathbb Z_2$ which satisfies $\A \prec_\M \L(\Ga_{\sG_0})$. Passing to the relative commutants intertwining we get $\L(\Ga_{{\rm lk}(\sG_0) })\subseteq \L(\Ga_{\sG_0})'\cap \M\prec_\M \A'\cap \M =\A$. Thus a corner of $\L(\Ga_{{\rm lk}(\sG_0) })$ is amenable which implies that $\Ga_{{\rm lk}(\sG_0) }$ is amenable, a contradiction. $\hfill\qedsymbol$

\vskip 0.03in
\noindent\emph{Proof of Corollary \ref{nocartan2}.} As $\sG$ is not complete and $|\Ga_v|\geq 3$ for all $v\in \sV$ then $\sG$ does not have a de Rham join decomposition of the form  $\sG=\sG_1\circ \sG_2\circ \cdots \circ \sG_k$ so that for all $2\leq i\leq k$ we have $\sG_i=(\{v_i,w_i\}, \emptyset)$ and $\Ga_{v_i}=\Ga_{w_i}\cong \mathbb Z_2$. The result follows from Corollary \ref{nocartan1}. $\hfill\qedsymbol$

\vskip 0.03in
\noindent \emph{Proof of Theorem \ref{uniquecartan1}} By Theorem \ref{controlmasa2}, $\A\prec_\M L^\infty(X)$ and the conclusion follows from Theorem \ref{Po01}.
$\hfill\qedsymbol$

\vskip 0.05in
\noindent {\bf Remarks}. Both graph conditions presented in the previous corollary are optimal. In other words, for these graphs there there are choices of the vertex groups for which the result does not hold as stated. 

To see this, let first $\mathscr G= \{w\} \circ {\rm lk}(w)$ be a star-shaped graph.  Using a result from \cite{OP07}, we show there is a graph product group $G= \mathscr G\{\Ga_v\}= \Ga_w \times \Ga_{{\rm lk}(w)}$  and a pmp action $G \ca X$ such that  $\M=L^\infty(X)\rtimes G$ has at least two unitary non-conjugate Cartan subalgebras. Indeed, let $\Ga_w= \mathbb Z^2\rtimes {\rm SL}_2(\mathbb Z)$ and let $\Ga_w\ca Y=\underset{\longleftarrow}{\lim} (Y_n)$ be the free ergodic profinite action as in \cite[Theorem D]{OP08}. Also let  $\Ga_{{\rm lk}(w)}\ca Z$ be any free ergodic pmp action. Then the consider the product action $G=\Ga_w \times \Ga_{{\rm lk}(w)}\ca Y\times Z=:X$. Then clearly $L^\infty(X)$ and $\L(\mathbb Z^2)\bar\otimes L^\infty(Z)$ are two Cartan subalgebras of $\M$ which are not unitarily conjugated. 

For the second situation let $\mathscr G= (\{w, v\}, \emptyset) \circ {\rm lk}(w,v)$ and assume $\G_v=\G_w\cong \mathbb Z_2$. Thus we have a product decomposition of graph product group
$G=\mathscr G\{\Ga_v\}= (\mathbb Z_2\ast \mathbb Z_2)\times \Ga_{{\rm lk}(w)} $. Since the infinite dihedral group  $\Lambda =\mathbb Z_2\ast \mathbb Z_2$ is ressidually finite, there is a free ergodic profinite pmp action $\Lambda \ca Y$. Let  $\Ga_{{\rm lk}(w)}\ca Z$ be any free ergodic pmp action and consider the product action $G=\Ga_w \times \Ga_{{\rm lk}(w)}\ca Y\times Z=:X$.   Next we view the infinite dihedral group as a semidirect product group $\Lambda= \mathbb Z \rtimes \mathbb Z_2=\langle a,s \,|\, s^2=1, sas=a^{-1}\rangle $. Then one can check that $L^\infty(X)$ and $\mathcal L(\langle a\rangle )\bar\otimes L^\infty(Z)$ are two Cartan subalgebras of $\M$ that are not unitarily conjugated. The only noncanonical part of this statement is to show that $\mathcal L(\langle a\rangle )\bar\otimes L^\infty(Z)\subset \M$ is Cartan subalgebra. To see this it suffices to show that  $\mathcal L(\langle a\rangle )\subset L^\infty(Y)\rtimes \La =:\N$ is a Cartan subalgebra. Throughout the remainder of the proof let $\Omega = \langle a \rangle$ and first we show $\L(\Omega)$ is a masa. Let $x\in \L(\Omega)'\cap \N$. Hence $xu_g=u_gx$ for all $g\in\Omega$ and using the Fourier decomposition $x=\sum_{h\in \La}  x_h u_h$, where $x_h \in L^\infty(Y)$, we  get

\begin{equation}\label{commutation2}x_{ghg^{-1}}=\sigma_g(x_h)\text{ for all }g\in\Omega,h\in\Lambda.\end{equation}

Fix $a^ts=h\in \La\setminus\Omega$, where $t\in \mathbb Z$. Thus for every $g= a^r$, using the presentation of $\Lambda$ we have  $ghg^{-1}= a^{r+t}s a^{-r}= a^{r+t-1} sa^{-r-1}=\cdots = s a^{-2r-t}$. In particular this implies that $\{ghg^{-1} \,|\, g\in\Omega$ is infinite. Using \eqref{commutation2} we get that $x_h=0$ for all $h\in \La\setminus\Omega$. Now fix $h\in \langle a\rangle$. Since $\Omega$  is abelian then \eqref{commutation2} implies that 
$x_{h}=\sigma_g(x_h)$ and since $\sigma$ is free we conclude that $x_h\in \mathbb C 1$ for all  $h\in \Omega$. Altogether, these show that $x=\sum_h x_h u_h\in \L(\Omega )$ which give that $\mathcal L(\Omega )\subset \N$  is  a masa. Now we briefly argue it is a Cartan subalgebra. Since $\La \ca Y:=\underset{\longleftarrow}{\lim} (Y_n) $ is profinite then we have an increasing union of finite dimensional, abelian, $\Lambda$-invariant von Neumann algebras $\cdots \subset \B_n:=L^\infty(Y_n)\subset \B_{n+1}:=L^\infty(Y_{n+1})\subset \cdots \subset L^\infty(Y)$ with $\overline{\bigcup_n \B_n}^{sot}=L^\infty(Y)$. Let $\B_n =\oplus_{ 1\leq i\leq k_n} \mathbb C b^n_i$ for some $\tau$-orthonormal family $\{b^n_i\}\subset \B_n$ and observe that for every $h\in \Omega$ we have $\sigma_h(b^n_i)= \sum_j\mu^n_{i,j}(h) b^n_j$ where $\mu^n_{i,j}(h)=\langle \sigma_g(b^n_i), b^n_j\rangle\in \mathbb C$ with $|\mu^n_{i,j}(h)|\leq 1$. This clearly implies that   $L(\Omega) b_i^n\subseteq \sum^{k_n}_{j=1} b_j^n L(\Omega)$ and $b^n_iL(\Omega) \subseteq \sum^{k_n}_{j=1} L(\Omega) b_j^n$ for all $i$. In particular, for all $i$ and all $n$ we have $b_i^n\in \mathscr {QN}_\N(L(\Omega))$---the quasinormalizer of $L(\Omega)$ in $\N$, \cite[Section 1.4.2]{Po01}. Thus $\bigcup_n \B_n\subset \mathscr {QN}_\N(L(\Omega))$ and hence $L^\infty(Y)\subseteq \mathscr {QN}_\N(L(\Omega))''$. Since $\Omega$ is normal in $\Lambda$ we further get  $L^\infty(Y)\rtimes \Lambda\subseteq \mathscr {QN}_\N(L(\Omega))''$ and hence $\N= \mathscr {QN}_\N(L(\Omega))''$. However, since $L(\Omega)\subset \N$ is  masa, by \cite[Proposition 1.4.2]{Po01}, we have $\mathscr {QN}_\N(L(\Omega))''= \mathscr {N}_\N (\L(\Omega))''$ and hence $\L(\Omega)\subset \N$ is a Cartan subalgebra as claimed.

If one considers the star-shaped graphs excluded in the statement of Theorem \ref{uniquecartan1} one can still get a uniqueness of Cartan subalgebra result if stronger conditions are imposed on the vertex groups. To properly introduce the result we recall the following definition. 

Following \cite{PV11} a group $\Gamma$ is called $\C$-superrigid if the following property holds. Given any $\Gamma \ca \P$ trace-preserving action and   any projection $p\in \P\rtimes \Gamma=\M$ then any masa $\A\subseteq p\M p$ with finite index normalizer  $[p\M p: \mathscr N_{p \M p}(\A)'']<\infty$ must satisfy $\A\prec_{\M} \P$. From \cite{PV11, PV12} it follows that all free groups with at least two generators and, more generally, any non-elementery hyperbolic groups are $\C$-superrigid.

\begin{cor}\label{uniquecartan2} Let $\Gamma=\mathscr G\{\Gamma_v\}$ be any graph product where the vertex groups $\Gamma_v$ are non-amenable $\C$-superrigid, for all $v$. Then $\Gamma$ is a $\C$-superrigid group.  

\end{cor}

\begin{proof} Letting $\mathscr G=(\mathscr V,\mathscr E)$, we proceed by induction on $|\mathscr V|$. When $|\mathscr V|=1$ the conclusion follows from the definition. Thus it only remains to show the inductive step. If $\mathscr G$ is not star-shaped the result already follows from Theorem \ref{controlmasa2}. Therefore, since the vertex groups are infinite, we only have to treat the case $\mathscr G= \{v\}\circ {\rm lk}(v)$.  
Let $\Gamma \ca \P$ be any trace preserving action and let $0\neq p\in \P\rtimes \Gamma=:\M$ be any projection.  Let $\A\subseteq p\M p$ be a masa so that  $[p\M p: \mathscr N_{p \M p}(\A)'']<\infty$. Since $\Gamma = \Gamma_v \times \Gamma_{{\rm lk}(v)}$ we can view $\M = (\P \rtimes \Gamma_{{\rm lk}(v)})\rtimes \Gamma_v$ and since $\Gamma_v$ is $\C$-superrigid we further get that $\A \prec_\M \P \rtimes \Gamma_{{\rm lk}(v)}=:\N$.  As $\A\subset p\M p$ is a masa, using  \cite[Proposition 3.6]{CIK13} there are projections $0\neq a\in \A,0\neq e\in \N$, a masa $\B\subset e\N e$, a projection $0\neq e'\in \B'\cap e\M e$ and a unitary $u\in \M $ such that the inclusion $\Q:=\sN_{e \N e}(\B)''\subseteq e \N e$ has finite index, and the support satisfies $s(E_\N(e'))=e$ and $u\A a u^*=\B e'$. 
As $|\rm lk(v)|= |\mathscr V|-1$, the inductive hypothesis implies that $\Gamma_{{\rm lk}(v)}$  is $\C$-superrigid. Thus the prior relation further implies that $\B \prec_\N \P$.   Hence one can find a partial isometry $0\neq w\in \N$ with $w^*w\in \B$ such that $w\B \subseteq \P w$. Same argument from the proof of Theorem \ref{controlmasa1} shows that $we'\neq 0$.  Hence $w e'\B =w \B e'\subseteq \P we'$ and combining this with the previous containment it yields $w u a\A \subseteq \P we'u$. In particular, $\A \prec_\M \P$ which concludes the proof of the inductive step. \end{proof}

\end{document}